# EXPÉRIMENTATION D'UNE RESSOURCE POUR UNE SITUATION DE RECHERCHE ET DE PREUVE ENTRE PAIRS

Jean-Philippe GEORGET[*] – Baptiste LABROUSSE[**]

**Résumé –** Cette contribution présente une expérimentation testant une ressource destinée à un enseignant de l'école primaire. Elle fait suite à une recherche présentée à EMF2009. La ressource doit permettre la mise en œuvre d'une situation de recherche et de preuve entre pairs par un enseignant novice dans la mise en œuvre de ce type de situations. La contribution aborde les problèmes posés par ce projet puis explicite les moyens mis en œuvre pour élaborer la ressource et la méthodologie d'expérimentation dans une classe. Enfin, des résultats sont donnés en illustrant la complexité des processus en jeu.

**Mots-clefs** : ressource destinée aux enseignants, ergonomie, situations de recherche et de preuve entre pairs (RPP)

**Abstract** – This paper is about an experiment with the goal of testing a primary teacher's resource. This is the next part of a research presented at EMF2009. This resource must help a teacher to practice a research and proof activity between peers with her pupils even she had never done this before. The paper presents the problems given by this project and describes the means to build the resource and the methodology of the experiment with some classmates. The results are given illustrating the complexity of the underlying process.

Depuis de nombreuses années, les programmes scolaires de nombreux pays, la France en particulier, demandent aux enseignants de faire vivre à leurs élèves des situations de recherche et de preuve entre pairs (situations RPP) en classe de mathématiques. Pour autant, malgré les demandes institutionnelles renouvelées, les formations d'enseignants mises en place, les recherches et les expérimentations scientifiques menées dans le domaine, les pratiques ordinaires n'évoluent pas sensiblement (Artigue et Houdement 2007, Georget 2009, 2010). Cet état de fait témoigne en partie de la complexité du projet consistant à provoquer des évolutions des pratiques enseignantes en ce domaine. Deux éléments principaux contribuent à l'expliquer, la complexité des situations RPP et le manque d'ergonomie des ressources enseignantes (Georget 2009, 2010).

Une première expérimentation s'est attachée à étudier sur trois ans des moyens de favoriser la pratique de situations RPP par une dizaine d'enseignants de l'école primaire quasiment tous novices dans cette pratique (Georget 2009). Il s'agissait de s'appuyer sur l'émergence d'une communauté de pratique (Wenger 1998). L'expérimentation présentée ici (menée par un étudiant de Master d'enseignement) s'appuie sur ces travaux exploratoires et s'attache à étudier comment une ressource peut aider un enseignant à mettre en œuvre une situation de recherche et de preuve entre pairs donnée. Le contexte contraint de ce travail explique en grande partie le fait que la méthodologie soit relativement restreinte par endroits.

## I.    LES SITUATIONS RPP

L'expression *situations de recherche et de preuve entre pairs* (situations RPP) a été proposée pour désigner des situations dans lesquelles on place les élèves ou les étudiants dans des situations proches de celle d'un mathématicien lorsqu'il cherche un problème nouveau (Georget 2009). Cette nouvelle expression regroupe diverses situations de classe comme les *problèmes ouverts* (Arsac et Mante 2007), les *situations de recherche en classe* (Grenier et Payan 2002) et d'autres encore. Ces dernières appellations désignent peu explicitement les objectifs visés par ces situations par lesquelles on souhaite développer chez les élèves, voire parfois chez les enseignants, des savoirs et des savoir-faire propres aux débats mathématiques (par exemple le statut des conjectures, des preuves, des exemples et contre-exemples, etc.). En particulier, la qualification RPP a pour premier intérêt d'évoquer explicitement la composante « preuve » de ces situations, caractéristique primordiale fréquemment absente des déroulements dans les classes ordinaires, ainsi que la composante sociale « entre pairs », fréquemment délaissée, elle aussi, au profit de preuves relativement formelles dont la logique et le déroulement échappent à nombre d'élèves. Le second intérêt de cette appellation est de désigner une classe de situations similaires — posant de ce fait la question de leur efficacité respective. Une étude comparative a été menée en se basant sur une revue de littérature (Georget 2009). Elle a mis en évidence la particulière pertinence de l'approche des *problèmes ouverts* proposée par Arsac et Mante (2007). En effet, celle-ci nécessite peu de moyens à mettre en œuvre pour espérer obtenir des résultats similaires à ceux des autres approches — modulo le fait que l'évaluation de l'efficacité des différentes approches restent particulièrement sujette à caution du fait des méthodologies employées ou du manque d'explicitation de ces méthodologies dans la littérature.

---

[*] Université de Caen Basse-Normandie, CERSE, EA 965 – France – jean-philippe.georget@unicaen.fr
[**] Université de Tours – France – baptiste.labrousse@etu.univ-tours.fr

Il ressort de la même étude que des situations RPP peuvent être proposées en classe dès l'enseignement primaire. Se pose alors la question du manque de diffusion de ces pratiques de classe hors des cadres expérimentaux. Les réponses sont multiples et cette contribution ne cherchera pas à en dresser une liste. Nous nous limiterons aux réponses envisageables concernant les caractéristiques de ces situations, leurs potentiels, et l'ergonomie des ressources destinées aux enseignants.

## .1   Potentiels des situations RPP et processus dynamiques en cours de séance

Les situations proposées aux élèves doivent bien sûr permettre aux élèves de chercher. Ceci suppose que la solution ne soit pas déjà connue des élèves concernés, qu'elle ne paraisse pas évidente et peut-être aussi que plusieurs pistes paraissent intéressantes à explorer aux yeux des élèves, même si une unique solution s'avère pertinente au final. C'est le potentiel de recherche de la situation (Georget 2009). Il paraît aussi nécessaire que la situation ne soit pas résolue trop rapidement par les élèves, c'est à dire qu'elle leur résiste. C'est son potentiel de résistance. Si le problème paraît toujours insoluble tout au long de la recherche des élèves, ces derniers risquent de se démobiliser. La résistance de la situation RPP face aux tentatives des élèves doit donc varier, ce que nous appelons son potentiel de résistance dynamique. Enfin, la situation doit permettre à l'élève d'apprendre quelque chose, sous peine de quoi sa présence à l'école deviendrait problématique. C'est son potentiel didactique. Toutes ces caractéristiques s'actualisent différemment durant un déroulement de séance donné en fonction de multiples facteurs et dans des processus dynamiques qui ne sont pas toujours contrôlables, d'où l'usage du terme *potentiel* (Georget 2009).

Les situations de recherche ne permettent pas toujours à l'enseignant de mener des débats mathématiques dans sa classe, c'est pourquoi nous avons aussi défini la notion de potentiel de débat. La possibilité d'organiser de tels débats est susceptible d'évoluer au cours d'une séance et dépend, elle aussi, de multiples facteurs. En premier lieu, pour assurer un potentiel de débat intéressant, le problème choisi peut par exemple permettre à plusieurs solutions différentes d'apparaître, ce qui renforcera la conviction chez les élèves que des argumentations puis des preuves ou des ébauches de preuve sont nécessaires. L'expérience de l'enseignant et sa gestion de la séance sont d'autres facteurs permettant à ce potentiel de s'exprimer. Les consignes employées par les enseignants au moment où tout semble propice pour que des débats mathématiques émergent sont souvent non pertinentes. Par exemple, les enseignants demandent fréquemment aux élèves de « dire ce qu'ils ont fait », c'est à dire de raconter leur processus de recherche, mais sans pour autant leur demander de se prononcer sur la validité de la démarche et des résultats et sans solliciter des échanges entre pairs. *In* fine c'est alors à l'enseignant d'initier et de gérer des processus de preuve par des questions fermées. Sans doute aussi fréquemment, les enseignants valident certaines productions d'élèves avant toute mise en commun des productions, parfois involontairement ou implicitement. Les élèves sont alors moins motivés pour défendre leur travail puisque, bien que comptant sur eux pour lancer des processus de preuve, l'enseignant se sera prononcé sur la validité de leur travail. L'expérience des élèves est aussi un élément important qui joue dans l'actualisation du potentiel de débat. Par exemple, même si les consignes de l'enseignant ne sont pas adéquates, les élèves qui ont l'expérience des situations RPP ne les suivront pas à la lettre et auront un comportement pertinent (Georget 2009). À l'inverse, si les élèves n'ont pas suffisamment d'expérience, ils ne pourront pas anticiper la démarche sous-jacente à la situation.

En conclusion de cette partie, il est donc évident que, pour favoriser la pratique des activités RPP dans les classes, il faut proposer aux enseignants des situations ayant les potentiels suffisants pour le faire. Sachant que l'expression de ces potentiels dépend des déroulements des séances, il faut aussi proposer ces situations sous la forme de ressources permettant aux enseignants d'avoir l'information suffisante pour les exploiter, c'est à dire ayant des caractéristiques ergonomiques suffisantes.

*2. Ergonomie des ressources*

Des concepts utilisés dans les recherches sur l'ergonomie des environnements informatiques pour l'apprentissage humain (EIAH) sont pertinents pour analyser et élaborer des ressources destinées aux enseignants, que ces ressources soient numériques ou non (Georget 2010). Les concepts d'utilité et d'utilisabilité sont deux d'entre eux. L'évaluation de l'utilité d'une ressource destinée à un enseignant consiste à évaluer si elle lui permet effectivement d'atteindre l'objectif que l'on a fixé a priori à la ressource, c'est à dire ici de lui permettre de faire vivre aux élèves une situation RPP. Quant à l'évaluation de l'utilisabilité d'une ressource, elle consiste à évaluer si l'enseignant peut l'adapter à sa pratique et si elle prend en compte son degré d'expertise. L'évaluation de l'utilité et de l'utilisabilité peut s'effectuer *a priori*, c'est à dire hors de toute utilisation réelle de la ressource, et *a posteriori*, c'est à dire par l'observation et l'analyse des effets produits par son utilisation réelle par un enseignant.

Une étude de ressources existantes a montré la possibilité d'améliorer leur utilité et leur utilisabilité en proposant, implicitement ou explicitement, des options d'utilisation aux enseignants là où les ressources ne donnent souvent qu'une seule façon de mener une situation donnée (Georget 2009). C'est en effet un moyen possible de favoriser des genèses documentaires prolixes (Gueudet et Trouche 2008) et de provoquer des exploitations — optimales ou minimales — des situations mathématiques proposées, c'est à dire en préservant au moins leur statut de situation RPP. Ceci semble particulièrement pertinent lorsque l'on s'adresse à des enseignants expérimentés mais novices en matière de situations RPP, ceux-ci ne souhaitant pas toujours qu'une ressource leur dise ce qu'ils ont à faire de manière plus ou moins péremptoire — à tort ou à raison.

I. UNE NOUVELLE EXPÉRIMENTATION

Dans l'expérimentation menée précédemment (Georget 2009), un collectif d'enseignants était concerné. La présente expérimentation concerne un unique enseignant de CM1 (9-10 ans), volontaire, bénévole et expérimenté (10 ans de pratique). Il est relativement novice dans la pratique des situations RPP et il lui a juste été demandé de lire la ressource et de mettre en œuvre le problème comme il le souhaitait. Contrairement à la première expérimentation, il n'a pu s'appuyer ni sur un temps long ni sur une communauté de pratique initiée pour l'occasion pour faire évoluer sa pratique. Nous allons ici supposer en première approximation que la ressource proposée sera le principal élément d'aide pour lui permettre de mettre en place la situation RPP proposée. Dans ce contexte, la ressource préalablement utilisée dans la première expérimentation a été retravaillée.

   *1. Une ressource retravaillée*

Dans la précédente expérimentation, plusieurs informations étaient, d'une part, disponibles sur un site Web conçu à cet effet et, d'autre part, étaient discutées entre les participants durant les réunions qui ont eu lieu tout au long des trois années. Elles étaient aussi évoquées dans les comptes-rendus d'expérience qu'ont pu rédiger et échanger les enseignants. Pour résumer, ces informations étaient susceptibles de rester à l'esprit des enseignants pour plusieurs raisons et ont pu influencer leur pratique. Ici, le seul moyen de communiquer à l'enseignant des informations utiles est la ressource. Elle doit donc être claire, concise, précise, tout en offrant une utilisabilité supposée optimale (Georget 2010).

Cette ressource est basée sur le problème consistant à dénombrer les cordes joignant un nombre fixés de points disposés sur un cercle (ERMEL 1999). En plus de travailler des compétences liées aux situations RPP, les élèves peuvent développer des compétences relatives au dénombrement exhaustif d'objets vérifiant des propriétés données et à la généralisation à partir de quelques cas particuliers, toutes deux au cœur de cette situation. Nous avons vérifié que le problème offrait *a priori* des potentiels suffisants pour des élèves de la fin de l'école primaire (Georget 2009). Par ailleurs, le côté « abstrait » est intéressant, car il permet aux enseignants de constater que leurs élèves sont tout à fait capables de s'attaquer à ce type des situations sans qu'un contexte « concret » ne soit nécessaire, bien que cela reste envisageable.

   Pour l'essentiel, la structure de la ressource (cf. Annexe 1) est la même que celle de l'expérimentation précédente (Georget 2009, 2010). Elle est composée d'une description rapide destinée à l'enseignant, de la solution et des moyens de l'obtenir. Des éléments supplémentaires visent à informer l'enseignant mais ne lui donnent pas une façon unique de mener la séance. Par exemple, il est seulement suggéré de traiter le cas des six points puis celui de dix. L'enseignant a une liberté et les éléments d'information sont là, implicitement ou explicitement, pour le guider dans ses choix. Les éléments de débats possibles, par l'intitulé même de cette rubrique et par la formulation de ses items, laissent une large place à des adaptations. La précédente expérimentation ayant montré que les enseignants avaient largement apprécié nos choix, notamment celui d'une ressource rapide à consulter, ils ont été conservés.

   Certains éléments ont été ajoutés ou modifiés. L'énoncé destiné à l'enseignant a été complété par deux exemples d'énoncés destinés aux élèves. Le cas de quinze points est implicitement proposé en supplément des cas six et dix points. Il s'agit d'éléments améliorant l'utilité et l'utilisabilité de la ressource.

   La preuve par récurrence, initialement dans la rubrique *Autres éléments de l'activité*, a été déplacée pour homogénéiser la ressource.

   La nouvelle rubrique *Conseils pédagogiques* présente un modèle de déroulement de séance qui semble plus opérationnel que celui proposé par Arsac et Mante (2007). En particulier, il propose à l'enseignant de ne pas attendre trop longtemps avant de faire une première médiation pour s'assurer que les élèves ont bien compris le problème initial. En effet, nous avons mis en évidence que des enseignants ne clarifient pas toujours, volontairement ou non, certaines caractéristiques de la situation de départ, ce qui cause des dysfonctionnements parfois jusqu'à la fin des séances (Georget 2009). Ce serait le cas par exemple si l'enseignant ne clarifie pas le fait qu'un diamètre est une corde ou qu'un point peut être l'extrémité de plusieurs cordes. La ressource propose aussi à l'enseignant de ne pas valider de réponses avant la tenue d'un débat. Nous avons vu que cela peut être un élément important pour que des débats puissent naître et vivre dans la classe. Dans un registre plus pédagogique, il est aussi suggéré une utilisation des affiches plus raisonnable que celle souvent observée qui consiste à demander aux élèves d'écrire leur démarche, ceci étant généralement peu productif dans les faits.

## 2. Autres éléments méthodologiques

L'enregistrement vidéo a eu pour objectif d'enregistrer le son pendant le déroulement de la séance et de filmer le tableau car il s'avère parfois très difficile de noter ce qui est écrit dessus lors des phases de mises en commun. Bien qu'il accentue généralement le stress de l'enseignant, ce moyen permet une meilleure récolte de données qu'avec une simple prise de notes et permet ainsi une reconstruction plus fidèle du déroulement de la séance (Georget 2010).

L'analyse des données consiste en premier lieu à rédiger une narration de la séance (Roditi 2001) à partir des notes prises et de l'enregistrement, c'est à dire un texte racontant son déroulement qui donne une première description appréhendable de la séance. Nous repérons des épisodes (Robert 1999, Roditi 2001) qui correspondent à une activité observée ou potentielle des élèves et qui peuvent correspondre à plusieurs des tâches prescrites ou attendues par l'enseignant (Robert et Rogalski 2002). Nous résumons chaque épisode qui nous semble significatif et indiquons sa durée approximative dans une trame simplifiée de séance. En particulier, certains épisodes tels que des interventions relativement extérieures au déroulement de la résolution du problème ne sont pas retenues dans le corpus (personne frappant à la porte, récréation, etc.). Elle est simplifiée aussi dans le sens où certains épisodes sont compactés en phases afin de faciliter la description de la séance en termes de présentation des tâches prescrites (phases codées par Pres, cf. Tableau 1), de recherches autonomes de véritables problèmes mathématiques par les élèves (RechO), de mises en commun (MC) des productions d'élèves et de véritables échanges entre élèves, de « cours dialogués » (CD), c'est-à-dire d'échanges enseignant-élèves fortement dirigés par l'enseignant pour faire avancer la résolution du problème. Nous étudions ensuite l'évolution des différents potentiels de la situation RPP et le rôle effectif qui nous paraît être celui de l'enseignant dans ces évolutions. Nous en sommes à ce stade réduits à faire des hypothèses sur ses logiques d'action (Robert et Rogalski 2002) en établissant un rapport probable avec le contenu de la ressource. Nous étudions aussi la possibilité d'une optimisation de la gestion de la séance et de la ressource, toutes choses étant supposées égales par ailleurs.

## 3. Résultats

Dans cette section, la trame simplifiée de la séance est donnée dans le Tableau 1. Les durées des épisodes sont arrondies au quart de minute le plus proche. La séance est ensuite étudiée en lien avec la ressource proposée.

La trame simplifiée de la séance permet tout d'abord de constater qu'un certain nombre de conseils donnés dans la ressource semblent suivis. C'est le cas (phase 1) de la présentation orale suivi d'un texte écrit, pratique peu courante à ce niveau d'enseignement chez les enseignants novices (Georget 2009). Par contre, il faut noter que la présentation d'une corde faite par l'enseignant n'est pas optimale puisqu'il ne précise qu'après une douzaine de minutes qu'une corde se trace avec la règle. De plus, la distinction corde-diamètre n'est pas faite clairement. Tout ceci aboutit à des problèmes de compréhension de l'énoncé qui pourraient être évités et qui se retrouvent tout au long de la séance, jusqu'à l'exploitation collective des productions. À cet égard, la ressource ne remplit que partiellement son rôle de ressource utile.

| Phases | Codes | Durées |
|---|---|---|
| 1. Présentation orale du cas des 6 points (tâche T6), mise au point du vocabulaire, reformulations de la consigne, écriture de la consigne au tableau | Pres | 6' |
| 2. Recherche individuelle de T6 | RechO | 3'30" |
| 3. Reformulations de la consigne | Pres | 2' |
| 4. Continuation de la recherche de T6 | RechO | 30" |
| 5. Précision : la corde se trace avec une règle ! | Pres | 0" |
| 6. Continuation de la recherche de T6 | RechO | 1'45" |
| 7. Présentation du cas général et du cas 15 points (tâches Tn, T15) | Pres | 4'30" |
| 8. Recherche individuelle de Tn et T15 | RechO | 5'30" |
| 9. Présentation de la recherche de Tn et T15 en binômes | Pres | 2'30" |
| 10. Recherche de Tn et T15 en binômes | RechO | 15" |
| 11. Reformulations de la consigne | Pres | 45" |
| 12. Continuation de la recherche Tn et T15 en binômes | RechO | 1'30" |
| 13. Annonce de la présentation à venir des productions des binômes | Pres | 30" |
| 14. Continuation de la recherche Tn et T15 en binômes | RechO | 5'30" |
| 15. Présentation quasi exhaustive des productions des binômes (NB. alternances MC/CD non détaillées ici) | MC/CD 12'/8' | 20' |
| Durée totale approximative | | 54'45" |

*Tableau 1 – Trame simplifiée de la séance observée.*

Cependant, on trouve aussi (succession des phases 1, 2 et 3) la présence d'une reformulation de la consigne, précédée d'une première courte recherche des élèves et suivie d'une autre recherche autonome des élèves. Bien qu'elle semble relever du « bon sens », cette pratique est peu courante chez les enseignants (Georget 2009). Même si on ne peut le certifier et même si l'enseignant nous a signalé ses difficultés d'organisation dans ce type de séances, il est donc possible que l'enseignant ait adopté le modèle de déroulement de séance proposé dans la ressource. La trame permet ensuite de constater l'absence de tâche fermée lors des épisodes recherches autonomes. Les potentiels de la situation sont donc assez bien préservés jusqu'à l'avant-dernière phase (phase 14). En particulier, l'objectif de faire travailler les élèves sur une véritable situation de recherche est atteint.

À l'inverse, le potentiel de débat et le potentiel didactique ne s'expriment pas ou alors pas de manière optimale. En effet, l'élaboration de preuves entre pairs est limitée par un ou des épisodes de cours dialogué (CD) pour chaque production. L'étude détaillée de la phase d'exploitation des productions (phase 15, 20') conclut à une durée totale de 8' pour le code CD contre 12' pour le code MC, caractéristique de moments plus favorables à l'élaboration de preuves entre pairs. Il faut retenir ici que l'enseignant prend largement la parole pour diriger la séance, alors qu'il pourrait davantage favoriser des débats mathématiques dans la classe. L'enseignant déclare qu'il « parle trop », constat courant chez les enseignants novices, mais c'est ici la pertinence du contenu de ses interventions qui empêche le potentiel de débat de s'exprimer. Du fait de l'enseignant ou des élèves eux-mêmes, plusieurs occasions existent d'orienter le déroulement de la séance vers un fonctionnement plus en phase avec les objectifs des situations RPP. Ces occasions n'aboutissent pas du fait d'un manque de soutien de l'enseignant. À cet égard, la ressource se révèle peu utile et peu utilisable quand elle évoque un rôle « d'animateur » sans en donner les clés élémentaires. En particulier, l'optimisation de la présentation des productions pouvait se faire ici en tenant compte de résultats redondants et divergents au lieu de poursuivre une présentation exhaustive. Le choix d'avoir une ressource concise a bien sûr pesé sur ces manques, mais certains gagneraient visiblement à être comblés. L'autre facteur explicatif de la main mise de l'enseignant sur les débats est sa volonté manifeste de faire résoudre le problème en une unique séance. En témoigne par exemple l'introduction (phase 7) de la tâche Tn de recherche du nombre de cordes en fonction du nombre de points dans le cas général, alors que les élèves n'ont même pas encore débattu du cas des 6 points, débat potentiellement prometteur puisque tous les élèves n'avaient toujours pas le même résultat à l'issue de la séance. Le déroulement de la séance révèle aussi une tension tenace entre cette volonté de tenir les délais et la volonté de laisser les élèves débattre et avancer eux-mêmes dans la résolution du problème. Sur ce point, la façon de compléter la ressource de façon utile, utilisable et acceptable reste pour nous ouverte. Enfin, vraisemblablement à cause du temps disponible, aucune phase de conclusion n'est présente, ne permettant pas au potentiel didactique de s'exprimer clairement, même *a minima*. La ressource n'évoquait pas ce point.

## II.    CONCLUSION

Dans l'expérimentation décrite, il était attendu que la ressource soit utile, c'est à dire permette à un enseignant de mettre en œuvre une situation de recherche et de preuve entre pairs dans sa classe malgré son manque d'expérience en ce domaine. Des dynamiques problématiques ont subsisté. En particulier, si les potentiels de recherche, de résistance et de résistance dynamique se sont bien exprimés, le potentiel de débat et le potentiel didactique ont été bridés par des interventions de l'enseignant. Les élèves ont donc bien vécu une situation de recherche mais ils n'ont que peu ou pas été en situation d'élaborer des preuves entre pairs. Ceci témoigne, d'une part, de la complexité intrinsèque de ces séances (Georget 2009) — l'enseignant ayant été confronté à des dynamiques problématiques qu'il n'a pas su gérer de manière optimale —, et d'autre part, de la présence de défauts ergonomiques dans la ressource qui lui a été proposée — défauts dont certains peuvent probablement être comblés facilement sans remettre en cause l'acceptabilité de la ressource. Enfin, même si le cadre d'expérimentation a été contraint (contexte de travail d'un étudiant de Master d'enseignement), la présente contribution montre combien l'élaboration d'une ressource ergonomique peut être un processus complexe. Les outils pour élaborer ce type de ressources semblent pertinents mais restent à affiner pour obtenir des résultats qui ne se bornent pas à des hypothèses, parfois périlleuses, concernant l'impact d'une ressource donnée sur la pratique d'un enseignant donné. En effet, nous avons considéré une ressource comme facteur explicatif essentiel de l'activité de l'enseignant mais la méthodologie employée ne permet pas de prouver la pertinence de cette analyse, elle ne permet que de l'envisager.

ANNEXE 1

1. *Le problème*

*Présentation*

Énoncé pour l'enseignant : on place un certain nombre de points sur un cercle. Est-il possible de trouver le nombre de cordes (segment joignant deux points du cercle) ?

Énoncés possibles pour les élèves :

- on place 6 points sur un cercle. Est-il possible de trouver le nombre de cordes ?
- on place 15 points au hasard sur un cercle. Peut-on trouver à coup sûr le nombre de cordes ?

*Exemples*

On peut commencer par 6 points disposés de façon irrégulière sur un cercle. On obtient 15 cordes. On peut ensuite passer à 10 points ce qui donne 45 cordes. Les élèves ont peu de chance de pouvoir les compter de façon sûre.

*Solutions*

Si on place $n$ points sur un cercle, le nombre de cordes est égal à : $n(n-1)/2$. Par exemple, pour 6 points, le nombre de cordes est égal à $6*5/2= 15$.

2. *Les preuves*

*Preuve additive*

La preuve revient à calculer la somme $(n-1) + \cdots + 3 + 2 + 1$.

En effet, on choisit un des points. Il permet d'obtenir $(n-1)$ cordes. En prenant un autre point, on obtient une corde de moins, c'est à dire $(n-2)$ et ainsi de suite jusqu'à l'avant-dernier point qui ne peut être joint qu'au dernier point, ce qui donne une seule corde.

Le calcul de la somme s'effectue de la manière suivante :
$(n-1) + (n-2) + \cdots + 3 + 2 + 1$
$+ 1 + 2 + \cdots + (n-3) + (n-2) + (n-1) = n + n + \cdots + n + n + n = n(n-1)$

On a alors calculé la somme 2 fois, il faut donc diviser l'expression par 2 pour obtenir le résultat.

*Preuve multiplicative*

Il y a $n$ points. Chaque point est relié à $(n-1)$ points. Mais, avec cette méthode, chaque corde est comptée 2 fois (une fois par extrémité). On obtient donc $n(n-1)/2$ cordes.

*Preuve par récurrence*

Une variante (raisonnement par récurrence) consiste à déduire par exemple le cas des 6 points de celui des 5 points (que l'on peut calculer à part). Il suffit d'ajouter les 5 nouvelles cordes créées par le sixième point à celles déjà comptées pour les 5 premiers (au nombre de 10). Ainsi, on obtient 5+10 = 15 cordes.

Le fait de proposer successivement certains cas, par exemple le cas de 5 puis de 6 points, risque d'induire cette stratégie.

*Éléments de débats possibles*

- méthode pour être sûr de compter toutes les cordes sans en oublier ;
- moyen de communiquer sa démarche (les élèves peuvent proposer plusieurs types de codage) ;
- les élèves peuvent proposer la preuve basée sur la multiplication sans qu'ils soient capables de l'expliquer dans un premier temps : cette preuve ne peut pas être considérée comme valide sans une explication acceptée par la classe. À l'issue des débats, l'enseignant peut en proposer une ;
- efficacité des différentes formules : elle dépend du nombre de points considéré.

*Autres éléments de l'activité*

- certains élèves risquent de confondre cordes et diamètres ;
- un nombre élevé de points oblige la recherche d'une méthode générale ;
- proposer des cas qui se succèdent (5 points, 6 points, etc.) risque d'induire la preuve par récurrence.

3. Conseils pédagogiques

Les conseils suivants ne constituent pas une obligation mais une aide potentielle à des difficultés rencontrées lors des situations de recherche. Libre à l'utilisateur de les suivre ou pas.

La première recherche individuelle sert à trouver des méthodes en vue de la recherche par groupe où elles seront discutées. Une médiation est utile si les élèves n'ont pas compris l'énoncé et s'il y a des précisions à apporter au problème, d'où un retour à la recherche individuelle puisque l'énoncé a changé.

La séance peut prendre la forme suivante :

- présentation du problème ;
- recherche individuelle (relativement courte) ;
- médiation ;
- suite de la recherche individuelle ou nouvelle recherche individuelle ;
- recherche en groupe ;
- mise en commun.

Énoncer le problème à l'oral pour éviter la barrière de l'écrit, le problème pourra être noté par la suite au tableau pour que les élèves n'oublient pas ce qui est cherché.

Composer des groupes de 3-4 élèves maximum pour que les échanges soient productifs et tenter de limiter un rapport de dominants/dominés, dans un binôme, ainsi que le nombre de réponses à exploiter lors des mises en commun.

Ne pas donner des éléments de réponse pendant les phases de recherche (en validant ou en invalidant par la voix ou la gestuelle), ceci évite que certains soient avantagés par rapport à d'autres, et évite des soucis lors de la mise en commun.

Le but est de maintenir le suspens quant à la résolution : ne pas attendre qu'un groupe ou un élève ait trouvé une méthode fiable pour avancer dans les phases si cela induit la fin prématurée de la séance.

Si un élève pense avoir trouvé toutes les solutions (alors que ce n'est pas le cas), soit il est possible de lui dire qu'il en reste ce qui revient à invalider la réponse, soit se contenter de lui demander s'il est sûr et lui montrer en fin de séance que la prochaine fois, il faudra d'avantage chercher.

Adopter une attitude « d'animateur » pendant la phase collective, l'important étant que les élèves se convainquent entre eux.

Les élèves voudront généralement tout expliquer par écrit, ce qui est compréhensible dans ce genre de situations. Une affiche peut constituer une simple aide à la présentation orale des solutions (pour éviter de copier l'ensemble des solutions, l'essentiel étant la démarche).